\documentclass[11pt,leqno]{article}
\normalfont
\renewcommand{\baselinestretch}{1.25}
\usepackage{amsfonts}

\newfont{\eulercursive}{eurm10 at 12pt}

\newcommand{\mya}{\mbox{\eulercursive a}}
\newcommand{\myb}{\mbox{\eulercursive b}}
\newcommand{\myc}{\mbox{\eulercursive c}}
\newcommand{\myd}{\mbox{\eulercursive d}}

\newcommand{\selt}{\mathbf{s}}

\newcommand{\QED}{\raisebox{0.5mm}{\fbox{\rule{0mm}{1.5mm}\ }}}

\newcommand{\LemmaOne}{Lemma 3.1} 
\newcommand{\LemmaOnePrime}{Lemma 3.2} 
\newcommand{\BothLemmaOnes}{Lemmas 3.1 and 3.2}
\newcommand{\LemmaTwo}{Lemma 3.3}

\newcommand{\MainTheorem}{Theorem 3.4}

\newcommand{\LucasTheorem}{Theorem 4.1}
\newcommand{\LucasCorollary}{Corollary 4.2}
\newcommand{\OddCountingTheorem}{Theorem 4.3}
\newcommand{\CountingLemmaConverse}{Lemma 4.4}
\newcommand{\CountingLemmaCentral}{Lemma 4.5}
\newcommand{\CountingLemmaOddRow}{Lemma 4.6}
\newcommand{\CountingLemmaS}{Lemmas 4.5 and 4.6}
\newcommand{\CountingLemmaIdentity}{Lemma 4.7}
\newcommand{\CountingLemmaParity}{Lemma 4.8}

\begin{document}

\newpage
\setcounter{page}{1} 
\renewcommand{\baselinestretch}{1}

\vspace*{-0.7in}
\hfill {\footnotesize July 25, 2018}

\begin{center}
{\Large \bf Counting odd numbers in truncations of Pascal's triangle}

\vspace*{0.05in}
\renewcommand{\thefootnote}{1}
Robert G.\ Donnelly,\footnote{Department of Mathematics and Statistics, Murray State
University, Murray, KY 42071\\ 
\hspace*{0.25in}Email: {\tt rob.donnelly@murraystate.edu}} 
\renewcommand{\thefootnote}{2} 
\hspace*{-0.07in}Molly W.\ Dunkum,\footnote{Department of Mathematics, Western Kentucky University, Bowling Green, KY 42101\\ 
\hspace*{0.25in}Email: {\tt molly.dunkum@wku.edu}}
\renewcommand{\thefootnote}{3} 
\hspace*{-0.07in}Courtney George,\footnote{Department of Mathematics, University of Kentucky, Lexington, KY 40508\\ 
\hspace*{0.25in}Email: {\tt courtney.george@uky.edu}\ \ (Courtney is an erstwhile undergraduate student of Dr.\ Dunkum.)}
\renewcommand{\thefootnote}{4} 
\hspace*{-0.07in}and Stefan Schnake\footnote{Department of Mathematics, The University of Oklahoma, Norman, OK 73019\\ 
\hspace*{0.25in}Email: {\tt sschnake@ou.edu}\ \  (Stefan is an erstwhile undergraduate student of Dr.\ Donnelly.)}

\end{center}

\begin{abstract}
A ``truncation'' of Pascal's triangle is a triangular array of numbers that satisfies the usual Pascal recurrence but with a boundary condition that declares some terminal set of numbers along each row of the array to be zero.  
Presented here is a family of natural truncations of Pascal's triangle that generalize a kind of Catalan triangle.  
The numbers in each array are realized as differences of binomial coefficients; as counts of certain lattice paths and tableaux; and as entries of representing matrices for certain linear transformations of polynomial spaces.  
Lucas's theorem is applied to determine precisely those truncations for which the number of odd entries on each row is a power of two.  
\begin{center}


\end{center}
\end{abstract}

{\bf \S 1\ \ Introduction.} 
Observe that the following conditions ({\sl i}), ({\sl ii}), 
and ({\sl iii}) uniquely determine an integer-valued function $A$ on 
$\mathbb{Z} \times 
\mathbb{Z}$:  
({\sl i}) $A(0,0) = 1$, 
({\sl ii}) $A(n,k) = 0$ if $n<0$, 
$k<0$, or $k > \lfloor n/2 \rfloor$, 
and ({\sl iii}) $A(n,k) = A(n-1,k-1) + A(n-1,k)$ 
for all other integer pairs $(n,k)$ when $n>0$.  
The output numbers of interest are those within the triangular array $\mbox{\Large (}A(n,k)\mbox{\Large )}$ indexed by integer pairs $(n,k)$ for which $0 \leq k \leq n$. When we display this array, we get a kind of ``truncation'' of Pascal's triangle.  Here are the first ten rows: 
\begin{center}
\begin{tabular}{ccccccccccccccccccc}
 & & & & & & & & & 1 & & & & & & & & & \\
 & & & & & & & &1 & & 0 & & & & & & & & \\
 & & & & & & & 1 & & 1 & & 0 & & & & & & & \\
 & & & & & & 1 & & 2 & & 0 & & 0 & & & & & & \\
 & & & & &1 & & 3 & & 2 & & 0 & &0 & & & & & \\
 & & & & 1 & & 4 & & 5 & & 0 & & 0 & & 0 & & & & \\
 & & & 1 & &5 & & 9 & & 5 & & 0 & & 0 & & 0 & & & \\
 & & 1 & & 6 & & 14 & & 14 & & 0 & & 0 & & 0 & & 0 & & \\
 & 1 & & 7 & & 20 & & 28 & & 14 & & 0 & & 0 & & 0 & & 0 & \\
1 & & 8 & & 27 & & 48 & & 42 & & 0 & & 0 & & 0 & & 0 & & 0
\end{tabular}
\end{center} 
Note the appearance of the famous Catalan numbers in the two rightmost nonzero 
``columns'' of the array and as sums of the numbers on northwest-to-southeast diagonals.  
The nonzero entries of this array are called ``ballot numbers,'' as they count the number of ways one candidate can defeat another candidate in a two-person election, under certain constraints. 
For further explication of this and other well-known phenomena related to this Catalan array, see for example \cite{Aval} and references therein.  

In this paper we present a generalization of this array by simply and 
naturally varying the ``boundary condition'' ({\sl ii}) above.  
We will have one such array for each positive integer $t$, 
where $t$ identifies the first row of the array that no longer fully agrees with Pascal's triangle, i.e.\ the first ``truncated'' row.
So, for example, the $t=1$ array is the Catalan array depicted above.  
The nonzero numbers in these more general 
triangular arrays are shown to 
be differences of binomial coefficients as well as counts of certain 
lattice paths. 
The fourth author, in consultation with the first author, studied these arrays in an undergraduate student honors
thesis \cite{Schnake} concerning differential operators on function spaces.  
A version of 
the motivating problem of that thesis is presented below, and in \MainTheorem\ it is shown how 
the numbers in our truncated Pascal arrays are coefficients for certain polynomials 
which arise in the study of differential operators.  
However, in \cite{Reuveni}, Reuveni independently presented the so-called ``Catalan's trapezoids,'' which are the same as our truncated Pascal's triangles but indexed and formatted somewhat differently. 
In \cite{ReuveniEtAl}, these trapezoids are applied in a probabilistic analysis of certain lattice-gas flow models. 

We close this introduction with some descriptive comments about the content of the paper. 
We think these Pascal-like arrays are inherently pretty and provide for 
an excellent enumerative example or exercise: We have a recurrence, 
an explicit formula, combinatorial interpretations, and a polynomial 
algebra context for these numbers, as summarized in \MainTheorem. 
This would seem to place 
us well within an enumeratively salubrious environment as envisioned by 
Stanley in Chapter 1 of his classic text \cite{Stanley}. 
Our main result -- \OddCountingTheorem\ -- is a (modest) enumerative application of these arrays that generalizes the well-known problem of 
counting odd numbers in the rows of Pascal's triangle; this new theorem, which appears in \cite{George}, was obtained by the third author in consultation with the first and second authors. 
For Pascal's triangle, solutions to this odd-counting problem and other related problems are entertainingly recounted in \cite{Granville}. 
Our work in this paper leaves open the possibility of generalizing other such problems from Pascal's triangle to the truncations of Pascal's triangle presented here.

\vspace*{0.1in}
{\bf \S 2\ \  A family of truncations of Pascal's triangle.} 
For the rest of the paper, $t$ denotes a fixed 
positive integer, which we think of informally as designating the first truncated row of the associated Pascal-like integer array. 
Consider a function $\mya_{t}:\, \mathbb{Z} \times 
\mathbb{Z} \longrightarrow \mathbb{Z}$ satisfying ({\sl i}) $\mya_{t}(0,0) = 1$, 
({\sl ii}) $\mya_{t}(n,k) = 0$ if $n<0$, 
$k<0$, or $k > \min\left\{\lfloor \frac{n-1+t}{2} \rfloor , 
n\right\}$, 
and ({\sl iii}) $\mya_{t}(n,k) = \mya_{t}(n-1,k-1) + \mya_{t}(n-1,k)$ 
for all other integer pairs $(n,k)$.  Figure 2.1 displays the first ten rows of the array 
$\mbox{\Large (}\mya_{4}(n,k)\mbox{\Large )}$ when viewed as a  
truncation of Pascal's triangle. 
We call an array $\mbox{\Large (}\mya_{t}(n,k)\mbox{\Large )}_{0 \leq k \leq n < \infty}$ a {\em truncated Pascal's triangle}. Of course, the Catalan array $\mbox{\Large 
(}A(n,k)\mbox{\Large )}$ of \S 1
is just the $t=1$ version of $\mbox{\Large 
(}\mya_{t}(n,k)\mbox{\Large )}$. 
In the next section we offer various interpretations of and contexts for the numbers in these truncated Pascal's triangles.

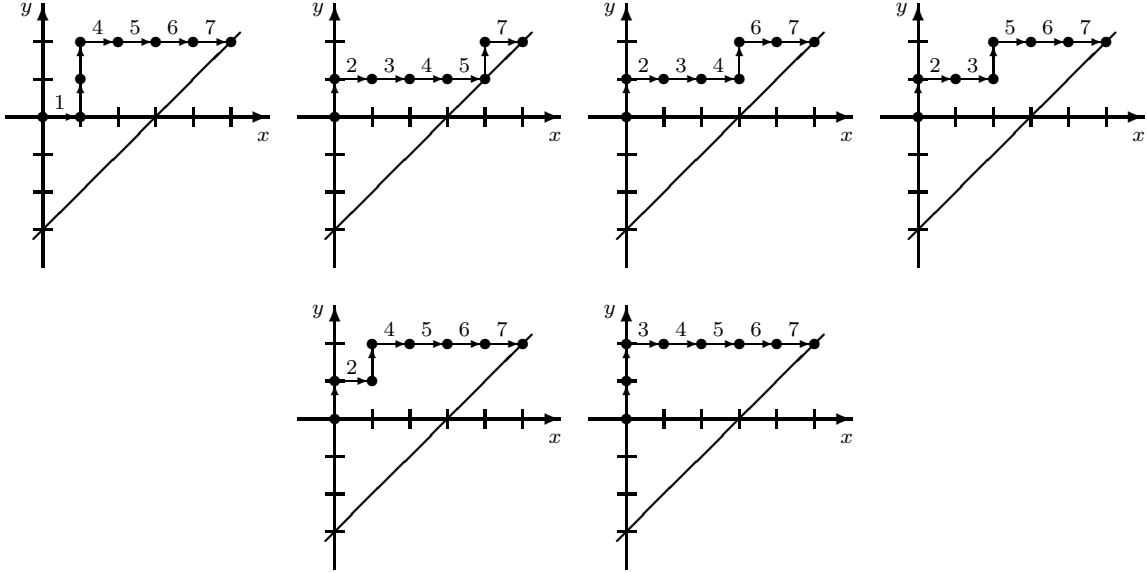
\begin{figure}[ht]
\begin{center}
\begin{tabular}{ccccccccccccccccccc}
 & & & & & & & & & 1 & & & & & & & & & \\
 & & & & & & & &1 & & 1 & & & & & & & & \\
 & & & & & & & 1 & & 2 & & 1 & & & & & & & \\
 & & & & & & 1 & & 3 & & 3 & & 1 & & & & & & \\
 & & & & &1 & & 4 & & 6 & & 4 & &0 & & & & & \\
 & & & & 1 & & 5 & & 10 & & 10 & & 4 & & 0 & & & & \\
 & & & 1 & &6 & & 15 & & 20 & & 14 & & 0 & & 0 & & & \\
 & & 1 & & 7 & & 21 & & 35 & & 34 & & 14 & & 0 & & 0 & & \\
 & 1 & & 8 & & 28 & & 56 & & 69 & & 48 & & 0 & & 0 & & 0 & \\
1 & & 9 & & 36 & & 84 & & 125 & & 117 & & 48 & & 0 & & 0 & & 0
\end{tabular}
\end{center} 

\centerline{Figure 2.1: The truncated Pascal's triangle $\mbox{\Large (}\mya_{4}(n,k)\mbox{\Large )}$.}
\end{figure}

\vspace*{-0.25in}
{\bf \S 3\ \  Algebraic-combinatorial aspects of truncated Pascal's triangles.} 
We aim to give several different descriptions of the numbers appearing in truncated Pascal's triangles, 
which will lead directly to \MainTheorem. 
Throughout this section, $n$ and $k$ are integers. 
Declare that  
\begin{equation}
\myb_{t}(n,k) := {n \choose k} - {n \choose k-t},
\end{equation}
with the usual understanding that the binomial coefficient ${p \choose q}$ is zero unless $0 \leq q \leq p$. 

Next, we count lattice paths.  An {\em NE-path} from $(x_{1},y_{1})$ 
to $(x_{2},y_{2})$ in the plane will be a continuous path starting 
at $(x_{1},y_{1})$, ending at $(x_{2},y_{2})$, and consisting of a 
finite number of unit steps in the north and east directions only. 
Say an NE-path from $(0,0)$ to $(k,n-k)$ is an $(n,k)$-NE-path, and call such a path {\em $t$-admissible}
if it does not intersect the line $y = x - t$; in such a case we say the path {\em stays weakly above} $y = x - t + 1$.  
For example, when $t=4$ and $(n,k)=(7,5)$, then $(k,n-k)=(5,2)$. 
As we can see in the pictures below, the number of $4$-admissible $(7,5)$-NE-paths is 14, i.e.\ there are fourteen NE-paths from $(0,0)$ to $(5,2)$ that stay weakly above $y=x-3$.  
For now, one can ignore the numbers assigned to the horizontal steps of each lattice path, although the pattern in which they are assigned should be evident. 
\begin{center}
\setlength{\unitlength}{0.5cm}
\begin{picture}(7.5,7)
\thicklines
\multiput(0.75,1)(0,1){6}{\line(1,0){0.5}}
\multiput(1,3.75)(1,0){6}{\line(0,1){0.5}}
\put(0.75,0.75){\line(1,1){5.5}}
\put(1,0){\line(0,1){7}}
\put(1,6){\vector(0,1){1}}
\put(0.4,6.8){\scriptsize $y$}
\put(0,4){\line(1,0){7}}
\put(6,4){\vector(1,0){1}}
\put(6.7,3.35){\scriptsize $x$}
\put(1,4){\circle*{0.3}}
\put(6,6){\circle*{0.3}}
\thinlines
\put(1,4){\vector(1,0){0.9}}
\put(2,4){\circle*{0.3}}
\put(2,4){\vector(1,0){0.9}}
\put(3,4){\circle*{0.3}}
\put(3,4){\vector(1,0){0.9}}
\put(4,4){\circle*{0.3}}
\put(4,4){\vector(0,1){0.9}}
\put(4,5){\circle*{0.3}}
\put(4,5){\vector(1,0){0.9}}
\put(5,5){\circle*{0.3}}
\put(5,5){\vector(0,1){0.9}}
\put(5,6){\circle*{0.3}}
\put(5,6){\vector(1,0){0.9}}
\put(1.3,4.2){\scriptsize $1$}
\put(2.3,4.2){\scriptsize $2$}
\put(3.3,4.2){\scriptsize $3$}
\put(4.3,5.2){\scriptsize $5$}
\put(5.3,6.2){\scriptsize $7$}
\end{picture}
\begin{picture}(7.5,7)
\thicklines
\multiput(0.75,1)(0,1){6}{\line(1,0){0.5}}
\multiput(1,3.75)(1,0){6}{\line(0,1){0.5}}
\put(0.75,0.75){\line(1,1){5.5}}
\put(1,0){\line(0,1){7}}
\put(1,6){\vector(0,1){1}}
\put(0.4,6.8){\scriptsize $y$}
\put(0,4){\line(1,0){7}}
\put(6,4){\vector(1,0){1}}
\put(6.7,3.35){\scriptsize $x$}
\put(1,4){\circle*{0.3}}
\put(6,6){\circle*{0.3}}
\thinlines
\put(1,4){\vector(1,0){0.9}}
\put(2,4){\circle*{0.3}}
\put(2,4){\vector(1,0){0.9}}
\put(3,4){\circle*{0.3}}
\put(3,4){\vector(1,0){0.9}}
\put(4,4){\circle*{0.3}}
\put(4,4){\vector(0,1){0.9}}
\put(4,5){\circle*{0.3}}
\put(4,5){\vector(0,1){0.9}}
\put(4,6){\circle*{0.3}}
\put(4,6){\vector(1,0){0.9}}
\put(5,6){\circle*{0.3}}
\put(5,6){\vector(1,0){0.9}}
\put(1.3,4.2){\scriptsize $1$}
\put(2.3,4.2){\scriptsize $2$}
\put(3.3,4.2){\scriptsize $3$}
\put(4.3,6.2){\scriptsize $6$}
\put(5.3,6.2){\scriptsize $7$}
\end{picture}
\begin{picture}(7.5,7)
\thicklines
\multiput(0.75,1)(0,1){6}{\line(1,0){0.5}}
\multiput(1,3.75)(1,0){6}{\line(0,1){0.5}}
\put(0.75,0.75){\line(1,1){5.5}}
\put(1,0){\line(0,1){7}}
\put(1,6){\vector(0,1){1}}
\put(0.4,6.8){\scriptsize $y$}
\put(0,4){\line(1,0){7}}
\put(6,4){\vector(1,0){1}}
\put(6.7,3.35){\scriptsize $x$}
\put(1,4){\circle*{0.3}}
\put(6,6){\circle*{0.3}}
\thinlines
\put(1,4){\vector(1,0){0.9}}
\put(2,4){\circle*{0.3}}
\put(2,4){\vector(1,0){0.9}}
\put(3,4){\circle*{0.3}}
\put(3,4){\vector(0,1){0.9}}
\put(3,5){\circle*{0.3}}
\put(3,5){\vector(1,0){0.9}}
\put(4,5){\circle*{0.3}}
\put(4,5){\vector(1,0){0.9}}
\put(5,5){\circle*{0.3}}
\put(5,5){\vector(0,1){0.9}}
\put(5,6){\circle*{0.3}}
\put(5,6){\vector(1,0){0.9}}
\put(1.3,4.2){\scriptsize $1$}
\put(2.3,4.2){\scriptsize $2$}
\put(3.3,5.2){\scriptsize $4$}
\put(4.3,5.2){\scriptsize $5$}
\put(5.3,6.2){\scriptsize $7$}
\end{picture}
\begin{picture}(7.5,7)
\thicklines
\multiput(0.75,1)(0,1){6}{\line(1,0){0.5}}
\multiput(1,3.75)(1,0){6}{\line(0,1){0.5}}
\put(0.75,0.75){\line(1,1){5.5}}
\put(1,0){\line(0,1){7}}
\put(1,6){\vector(0,1){1}}
\put(0.4,6.8){\scriptsize $y$}
\put(0,4){\line(1,0){7}}
\put(6,4){\vector(1,0){1}}
\put(6.7,3.35){\scriptsize $x$}
\put(1,4){\circle*{0.3}}
\put(6,6){\circle*{0.3}}
\thinlines
\put(1,4){\vector(1,0){0.9}}
\put(2,4){\circle*{0.3}}
\put(2,4){\vector(1,0){0.9}}
\put(3,4){\circle*{0.3}}
\put(3,4){\vector(0,1){0.9}}
\put(3,5){\circle*{0.3}}
\put(3,5){\vector(1,0){0.9}}
\put(4,5){\circle*{0.3}}
\put(4,5){\vector(0,1){0.9}}
\put(4,6){\circle*{0.3}}
\put(4,6){\vector(1,0){0.9}}
\put(5,6){\circle*{0.3}}
\put(5,6){\vector(1,0){0.9}}
\put(1.3,4.2){\scriptsize $1$}
\put(2.3,4.2){\scriptsize $2$}
\put(3.3,5.2){\scriptsize $4$}
\put(4.3,6.2){\scriptsize $6$}
\put(5.3,6.2){\scriptsize $7$}
\end{picture}
\end{center}

\begin{center}
\setlength{\unitlength}{0.5cm}
\begin{picture}(7.5,7)
\thicklines
\multiput(0.75,1)(0,1){6}{\line(1,0){0.5}}
\multiput(1,3.75)(1,0){6}{\line(0,1){0.5}}
\put(0.75,0.75){\line(1,1){5.5}}
\put(1,0){\line(0,1){7}}
\put(1,6){\vector(0,1){1}}
\put(0.4,6.8){\scriptsize $y$}
\put(0,4){\line(1,0){7}}
\put(6,4){\vector(1,0){1}}
\put(6.7,3.35){\scriptsize $x$}
\put(1,4){\circle*{0.3}}
\put(6,6){\circle*{0.3}}
\thinlines
\put(1,4){\vector(1,0){0.9}}
\put(2,4){\circle*{0.3}}
\put(2,4){\vector(1,0){0.9}}
\put(3,4){\circle*{0.3}}
\put(3,4){\vector(0,1){0.9}}
\put(3,5){\circle*{0.3}}
\put(3,5){\vector(0,1){0.9}}
\put(3,6){\circle*{0.3}}
\put(3,6){\vector(1,0){0.9}}
\put(4,6){\circle*{0.3}}
\put(4,6){\vector(1,0){0.9}}
\put(5,6){\circle*{0.3}}
\put(5,6){\vector(1,0){0.9}}
\put(1.3,4.2){\scriptsize $1$}
\put(2.3,4.2){\scriptsize $2$}
\put(3.3,6.2){\scriptsize $5$}
\put(4.3,6.2){\scriptsize $6$}
\put(5.3,6.2){\scriptsize $7$}
\end{picture}
\begin{picture}(7.5,7)
\thicklines
\multiput(0.75,1)(0,1){6}{\line(1,0){0.5}}
\multiput(1,3.75)(1,0){6}{\line(0,1){0.5}}
\put(0.75,0.75){\line(1,1){5.5}}
\put(1,0){\line(0,1){7}}
\put(1,6){\vector(0,1){1}}
\put(0.4,6.8){\scriptsize $y$}
\put(0,4){\line(1,0){7}}
\put(6,4){\vector(1,0){1}}
\put(6.7,3.35){\scriptsize $x$}
\put(1,4){\circle*{0.3}}
\put(6,6){\circle*{0.3}}
\thinlines
\put(1,4){\vector(1,0){0.9}}
\put(2,4){\circle*{0.3}}
\put(2,4){\vector(0,1){0.9}}
\put(2,5){\circle*{0.3}}
\put(2,5){\vector(1,0){0.9}}
\put(3,5){\circle*{0.3}}
\put(3,5){\vector(1,0){0.9}}
\put(4,5){\circle*{0.3}}
\put(4,5){\vector(1,0){0.9}}
\put(5,5){\circle*{0.3}}
\put(5,5){\vector(0,1){0.9}}
\put(5,6){\circle*{0.3}}
\put(5,6){\vector(1,0){0.9}}
\put(1.3,4.2){\scriptsize $1$}
\put(2.3,5.2){\scriptsize $3$}
\put(3.3,5.2){\scriptsize $4$}
\put(4.3,5.2){\scriptsize $5$}
\put(5.3,6.2){\scriptsize $7$}
\end{picture}
\begin{picture}(7.5,7)
\thicklines
\multiput(0.75,1)(0,1){6}{\line(1,0){0.5}}
\multiput(1,3.75)(1,0){6}{\line(0,1){0.5}}
\put(0.75,0.75){\line(1,1){5.5}}
\put(1,0){\line(0,1){7}}
\put(1,6){\vector(0,1){1}}
\put(0.4,6.8){\scriptsize $y$}
\put(0,4){\line(1,0){7}}
\put(6,4){\vector(1,0){1}}
\put(6.7,3.35){\scriptsize $x$}
\put(1,4){\circle*{0.3}}
\put(6,6){\circle*{0.3}}
\thinlines
\put(1,4){\vector(1,0){0.9}}
\put(2,4){\circle*{0.3}}
\put(2,4){\vector(0,1){0.9}}
\put(2,5){\circle*{0.3}}
\put(2,5){\vector(1,0){0.9}}
\put(3,5){\circle*{0.3}}
\put(3,5){\vector(1,0){0.9}}
\put(4,5){\circle*{0.3}}
\put(4,5){\vector(0,1){0.9}}
\put(4,6){\circle*{0.3}}
\put(4,6){\vector(1,0){0.9}}
\put(5,6){\circle*{0.3}}
\put(5,6){\vector(1,0){0.9}}
\put(1.3,4.2){\scriptsize $1$}
\put(2.3,5.2){\scriptsize $3$}
\put(3.3,5.2){\scriptsize $4$}
\put(4.3,6.2){\scriptsize $6$}
\put(5.3,6.2){\scriptsize $7$}
\end{picture}
\begin{picture}(7.5,7)
\thicklines
\multiput(0.75,1)(0,1){6}{\line(1,0){0.5}}
\multiput(1,3.75)(1,0){6}{\line(0,1){0.5}}
\put(0.75,0.75){\line(1,1){5.5}}
\put(1,0){\line(0,1){7}}
\put(1,6){\vector(0,1){1}}
\put(0.4,6.8){\scriptsize $y$}
\put(0,4){\line(1,0){7}}
\put(6,4){\vector(1,0){1}}
\put(6.7,3.35){\scriptsize $x$}
\put(1,4){\circle*{0.3}}
\put(6,6){\circle*{0.3}}
\thinlines
\put(1,4){\vector(1,0){0.9}}
\put(2,4){\circle*{0.3}}
\put(2,4){\vector(0,1){0.9}}
\put(2,5){\circle*{0.3}}
\put(2,5){\vector(1,0){0.9}}
\put(3,5){\circle*{0.3}}
\put(3,5){\vector(0,1){0.9}}
\put(3,6){\circle*{0.3}}
\put(3,6){\vector(1,0){0.9}}
\put(4,6){\circle*{0.3}}
\put(4,6){\vector(1,0){0.9}}
\put(5,6){\circle*{0.3}}
\put(5,6){\vector(1,0){0.9}}
\put(1.3,4.2){\scriptsize $1$}
\put(2.3,5.2){\scriptsize $3$}
\put(3.3,6.2){\scriptsize $5$}
\put(4.3,6.2){\scriptsize $6$}
\put(5.3,6.2){\scriptsize $7$}
\end{picture}
\end{center}

\begin{center}
\setlength{\unitlength}{0.5cm}
\begin{picture}(7.5,7)
\thicklines
\multiput(0.75,1)(0,1){6}{\line(1,0){0.5}}
\multiput(1,3.75)(1,0){6}{\line(0,1){0.5}}
\put(0.75,0.75){\line(1,1){5.5}}
\put(1,0){\line(0,1){7}}
\put(1,6){\vector(0,1){1}}
\put(0.4,6.8){\scriptsize $y$}
\put(0,4){\line(1,0){7}}
\put(6,4){\vector(1,0){1}}
\put(6.7,3.35){\scriptsize $x$}
\put(1,4){\circle*{0.3}}
\put(6,6){\circle*{0.3}}
\thinlines
\put(1,4){\vector(1,0){0.9}}
\put(2,4){\circle*{0.3}}
\put(2,4){\vector(0,1){0.9}}
\put(2,5){\circle*{0.3}}
\put(2,5){\vector(0,1){0.9}}
\put(2,6){\circle*{0.3}}
\put(2,6){\vector(1,0){0.9}}
\put(3,6){\circle*{0.3}}
\put(3,6){\vector(1,0){0.9}}
\put(4,6){\circle*{0.3}}
\put(4,6){\vector(1,0){0.9}}
\put(5,6){\circle*{0.3}}
\put(5,6){\vector(1,0){0.9}}
\put(1.3,4.2){\scriptsize $1$}
\put(2.3,6.2){\scriptsize $4$}
\put(3.3,6.2){\scriptsize $5$}
\put(4.3,6.2){\scriptsize $6$}
\put(5.3,6.2){\scriptsize $7$}
\end{picture}
\begin{picture}(7.5,7)
\thicklines
\multiput(0.75,1)(0,1){6}{\line(1,0){0.5}}
\multiput(1,3.75)(1,0){6}{\line(0,1){0.5}}
\put(0.75,0.75){\line(1,1){5.5}}
\put(1,0){\line(0,1){7}}
\put(1,6){\vector(0,1){1}}
\put(0.4,6.8){\scriptsize $y$}
\put(0,4){\line(1,0){7}}
\put(6,4){\vector(1,0){1}}
\put(6.7,3.35){\scriptsize $x$}
\put(1,4){\circle*{0.3}}
\put(6,6){\circle*{0.3}}
\thinlines
\put(1,4){\vector(0,1){0.9}}
\put(1,5){\circle*{0.3}}
\put(1,5){\vector(1,0){0.9}}
\put(2,5){\circle*{0.3}}
\put(2,5){\vector(1,0){0.9}}
\put(3,5){\circle*{0.3}}
\put(3,5){\vector(1,0){0.9}}
\put(4,5){\circle*{0.3}}
\put(4,5){\vector(1,0){0.9}}
\put(5,5){\circle*{0.3}}
\put(5,5){\vector(0,1){0.9}}
\put(5,6){\circle*{0.3}}
\put(5,6){\vector(1,0){0.9}}
\put(1.3,5.2){\scriptsize $2$}
\put(2.3,5.2){\scriptsize $3$}
\put(3.3,5.2){\scriptsize $4$}
\put(4.3,5.2){\scriptsize $5$}
\put(5.3,6.2){\scriptsize $7$}
\end{picture}
\begin{picture}(7.5,7)
\thicklines
\multiput(0.75,1)(0,1){6}{\line(1,0){0.5}}
\multiput(1,3.75)(1,0){6}{\line(0,1){0.5}}
\put(0.75,0.75){\line(1,1){5.5}}
\put(1,0){\line(0,1){7}}
\put(1,6){\vector(0,1){1}}
\put(0.4,6.8){\scriptsize $y$}
\put(0,4){\line(1,0){7}}
\put(6,4){\vector(1,0){1}}
\put(6.7,3.35){\scriptsize $x$}
\put(1,4){\circle*{0.3}}
\put(6,6){\circle*{0.3}}
\thinlines
\put(1,4){\vector(0,1){0.9}}
\put(1,5){\circle*{0.3}}
\put(1,5){\vector(1,0){0.9}}
\put(2,5){\circle*{0.3}}
\put(2,5){\vector(1,0){0.9}}
\put(3,5){\circle*{0.3}}
\put(3,5){\vector(1,0){0.9}}
\put(4,5){\circle*{0.3}}
\put(4,5){\vector(0,1){0.9}}
\put(4,6){\circle*{0.3}}
\put(4,6){\vector(1,0){0.9}}
\put(5,6){\circle*{0.3}}
\put(5,6){\vector(1,0){0.9}}
\put(1.3,5.2){\scriptsize $2$}
\put(2.3,5.2){\scriptsize $3$}
\put(3.3,5.2){\scriptsize $4$}
\put(4.3,6.2){\scriptsize $6$}
\put(5.3,6.2){\scriptsize $7$}
\end{picture}
\begin{picture}(7.5,7)
\thicklines
\multiput(0.75,1)(0,1){6}{\line(1,0){0.5}}
\multiput(1,3.75)(1,0){6}{\line(0,1){0.5}}
\put(0.75,0.75){\line(1,1){5.5}}
\put(1,0){\line(0,1){7}}
\put(1,6){\vector(0,1){1}}
\put(0.4,6.8){\scriptsize $y$}
\put(0,4){\line(1,0){7}}
\put(6,4){\vector(1,0){1}}
\put(6.7,3.35){\scriptsize $x$}
\put(1,4){\circle*{0.3}}
\put(6,6){\circle*{0.3}}
\thinlines
\put(1,4){\vector(0,1){0.9}}
\put(1,5){\circle*{0.3}}
\put(1,5){\vector(1,0){0.9}}
\put(2,5){\circle*{0.3}}
\put(2,5){\vector(1,0){0.9}}
\put(3,5){\circle*{0.3}}
\put(3,5){\vector(0,1){0.9}}
\put(3,6){\circle*{0.3}}
\put(3,6){\vector(1,0){0.9}}
\put(4,6){\circle*{0.3}}
\put(4,6){\vector(1,0){0.9}}
\put(5,6){\circle*{0.3}}
\put(5,6){\vector(1,0){0.9}}
\put(1.3,5.2){\scriptsize $2$}
\put(2.3,5.2){\scriptsize $3$}
\put(3.3,6.2){\scriptsize $5$}
\put(4.3,6.2){\scriptsize $6$}
\put(5.3,6.2){\scriptsize $7$}
\end{picture}
\end{center}

\begin{center}
\setlength{\unitlength}{0.5cm}
\begin{picture}(7.5,7)
\thicklines
\multiput(0.75,1)(0,1){6}{\line(1,0){0.5}}
\multiput(1,3.75)(1,0){6}{\line(0,1){0.5}}
\put(0.75,0.75){\line(1,1){5.5}}
\put(1,0){\line(0,1){7}}
\put(1,6){\vector(0,1){1}}
\put(0.4,6.8){\scriptsize $y$}
\put(0,4){\line(1,0){7}}
\put(6,4){\vector(1,0){1}}
\put(6.7,3.35){\scriptsize $x$}
\put(1,4){\circle*{0.3}}
\put(6,6){\circle*{0.3}}
\thinlines
\put(1,4){\vector(0,1){0.9}}
\put(1,5){\circle*{0.3}}
\put(1,5){\vector(1,0){0.9}}
\put(2,5){\circle*{0.3}}
\put(2,5){\vector(0,1){0.9}}
\put(2,6){\circle*{0.3}}
\put(2,6){\vector(1,0){0.9}}
\put(3,6){\circle*{0.3}}
\put(3,6){\vector(1,0){0.9}}
\put(4,6){\circle*{0.3}}
\put(4,6){\vector(1,0){0.9}}
\put(5,6){\circle*{0.3}}
\put(5,6){\vector(1,0){0.9}}
\put(1.3,5.2){\scriptsize $2$}
\put(2.3,6.2){\scriptsize $4$}
\put(3.3,6.2){\scriptsize $5$}
\put(4.3,6.2){\scriptsize $6$}
\put(5.3,6.2){\scriptsize $7$}
\end{picture}
\begin{picture}(7.5,7)
\thicklines
\multiput(0.75,1)(0,1){6}{\line(1,0){0.5}}
\multiput(1,3.75)(1,0){6}{\line(0,1){0.5}}
\put(0.75,0.75){\line(1,1){5.5}}
\put(1,0){\line(0,1){7}}
\put(1,6){\vector(0,1){1}}
\put(0.4,6.8){\scriptsize $y$}
\put(0,4){\line(1,0){7}}
\put(6,4){\vector(1,0){1}}
\put(6.7,3.35){\scriptsize $x$}
\put(1,4){\circle*{0.3}}
\put(6,6){\circle*{0.3}}
\thinlines
\put(1,4){\vector(0,1){0.9}}
\put(1,5){\circle*{0.3}}
\put(1,5){\vector(0,1){0.9}}
\put(1,6){\circle*{0.3}}
\put(1,6){\vector(1,0){0.9}}
\put(2,6){\circle*{0.3}}
\put(2,6){\vector(1,0){0.9}}
\put(3,6){\circle*{0.3}}
\put(3,6){\vector(1,0){0.9}}
\put(4,6){\circle*{0.3}}
\put(4,6){\vector(1,0){0.9}}
\put(5,6){\circle*{0.3}}
\put(5,6){\vector(1,0){0.9}}
\put(1.3,6.2){\scriptsize $3$}
\put(2.3,6.2){\scriptsize $4$}
\put(3.3,6.2){\scriptsize $5$}
\put(4.3,6.2){\scriptsize $6$}
\put(5.3,6.2){\scriptsize $7$}
\end{picture}
\end{center}
This count agrees with the $(n,k)=(7,5)$ entry of the example array $\mbox{\Large (}\mya_{4}(n,k)\mbox{\Large )}$ depicted in \S 2 above. 
Moreover, $\myb_{4}(7,5) = {7 \choose 5} - {7 \choose 1} = 14$. 
In general, we set 
\begin{equation}
\myc_{t}(n,k) := \mbox{\Large $|$}\{\mbox{\small  
$t$-admissible $(n,k)$-NE-paths}\}\mbox{\Large $|$}.
\end{equation}
For integers $x_{1}$, $y_{1}$, $x_{2}$, 
and $y_{2}$, the number of all NE-paths from $(x_{1},y_{1})$ to 
$(x_{2},y_{2})$ is easily seen to be $\displaystyle {y_{2}-y_{1} + 
x_{2}-x_{1} \choose x_{2}-x_{1}}$. 
Consider for the moment those NE-paths from $(0,0)$ to $(k,n-k)$ which 
cross the line $y = x - t + 1$ and therefore cross or touch the line $y = x - t$.  At the first point of intersection with $y = x - t$, reflect the initial part of the path across that line to obtain an NE-path from $(t,-t)$ to $(k,n-k)$.  This procedure can be reversed and therefore shows that the set of NE-paths from $(0,0)$ to $(k,n-k)$ which 
cross the line $y = x - t - 1$ is in one-to-one correspondence with the set of all 
NE-paths from $(t,-t)$ to $(k,n-k)$.  
(This is an instance of the famous reflection principle of Andr\'e.) 
That is, the number of NE-paths from $(0,0)$ to $(k,n-k)$ that cross the line $y=x-t-1$ is ${n-k+t+k-t \choose k-t}$.  
Therefore, 
\[\myc_{t}(n,k) = {n \choose k} - {n-k+t+k-t \choose k-t} = {n \choose k} - {n \choose k-t} = \myb_{t}(n,k).\] 
This shows: 

\noindent 
{\bf \LemmaOne}\ \ {\sl For all integers $n$ and $k$, we have $\myb_{t}(n,k) = \myc_{t}(n,k)$.}\hfill\QED 

For more about the well-known observation recorded above as \LemmaOne, see Chapter 1 of \cite{Mohanty}.  For a recent and readable survey of lattice path enumeration with references to many closely related results, see \cite{Humphreys}. 

In algebraic combinatorics, objects called ``tableaux'' are often used to index bases for representations of algebraic structures such as groups or Lie algebras.  
For tableaux relating to representations of symmetric groups, see for instance \cite{Sagan}; for many examples arising in Lie theory, see \cite{Proctor2}.   
Tableaux generally take the form of an array of boxes of some specified shape filled with integer entries subject to certain rules.  
Next, we offer a re-interpretation of $t$-admissible NE-paths as columnar tableaux. 
But before we do so, we very briefly remark on a connection between the Catalan array of \S 1 and the representation theory of symplectic Lie groups.  

The Lie theoretic notions of the present paragraph are mentioned in order to provide some further context for the kinds of objects we are considering, but readers unfamiliar with Lie representation theory can safely skip to the next paragraph. 
Taking $n \geq 1$, the number $A(2n+1,k) = \mya_{1}(2n+1,k)$ is the dimension of the $k$th fundamental representation of the symplectic group $\mbox{\em Sp}(2n,\mathbb{C})$, when $1 \leq k \leq n$.  
The number $A(2n,k) = \mya_{1}(2n,k)$ is the dimension of a certain indecomposable representation (that might reasonably be called the ``$k$th fundamental representation'') of the non-reductive odd symplectic group $\mbox{\em Sp}(2n-1, \mathbb{C})$; in the language of \cite{Proctor1}, this is the $\mbox{\em Sp}(2n-1,\mathbb{C})$-module corresponding to a partition of $k$ with $k$ parts. 
In fact, certain columnar tableaux provide a concrete connection between a number of the Catalan array (say, $A(2n,k)$ or $A(2n+1,k)$) and the dimension count of the associated symplectic group representation (say, the $k$th fundamental representation of $\mbox{\em Sp}(2n-1,\mathbb{C})$ or $\mbox{\em Sp}(2n,\mathbb{C})$): On the one hand, the tableaux index a special basis for the representing space, but on the other hand the number of these tableaux can easily be realized as the associated ballot number from the Catalan array. 

An $(n,k)${\em -columnar tableau} $T=(T_{1},T_{2},\ldots,T_{k})$ is a strictly increasing $k$-tuple of integers with $\{T_{1},T_{2},\ldots,T_{k}\} \subseteq \{1,2,\ldots,n\}$. 
We typically visualize such a tableau as a vertical column of $k$ boxes filled from top to bottom with the entries $T_{1},T_{2},\ldots,T_{k}$:
\begin{center}
\setlength{\unitlength}{0.5cm}
\begin{picture}(1,5)
\put(-1.8,2.6){$T = $}
\put(0,0){\line(0,1){5}} 
\put(1,0){\line(0,1){5}} 
\put(0,0){\line(1,0){1}} 
\put(0,1){\line(1,0){1}} 
\put(0,3){\line(1,0){1}} 
\put(0,4){\line(1,0){1}} 
\put(0,5){\line(1,0){1}} 
\put(0.2,4.3){\footnotesize $T_{1}$}
\put(0.2,3.3){\footnotesize $T_{2}$}
\put(0.375,2.5){\tiny $\bullet$}
\put(0.375,2.1){\tiny $\bullet$}
\put(0.375,1.7){\tiny $\bullet$}
\put(0.375,1.3){\tiny $\bullet$}
\put(0.2,0.3){\footnotesize $T_{k}$}
\end{picture}
\end{center}
An $(n,k)$-columnar tableau $T=(T_{1},\ldots,T_{k})$ is $t${\em -admissible} if, whenever $t-1 < k$, then $t-1+2j \leq T_{t-1+j}$ for any $j \in \{1,2,\ldots,k-t+1\}$. 
For example, when $t=4$ and $(n,k)=(7,5)$, then $k-t+1=2$, and our $j$'s are therefore from the set $\{1,2\}$. 
When $j=1$, $T_{t-1+j}=T_{4} \geq 5$, and when $j=2$, $T_{t-1+j}=T_{5} \geq 7$. 
So the $4$-admissible $(7,5)$-columnar tableaux are: 
\begin{center}
\setlength{\unitlength}{0.5cm}
\begin{picture}(1,5)
\put(0,0){\line(0,1){5}} 
\put(0.7,0){\line(0,1){5}} 
\put(0,0){\line(1,0){0.7}} 
\put(0,1){\line(1,0){0.7}} 
\put(0,2){\line(1,0){0.7}} 
\put(0,3){\line(1,0){0.7}} 
\put(0,4){\line(1,0){0.7}} 
\put(0,5){\line(1,0){0.7}} 
\put(0.2,4.3){\footnotesize $1$}
\put(0.2,3.3){\footnotesize $2$}
\put(0.2,2.3){\footnotesize $3$}
\put(0.2,1.3){\footnotesize $5$}
\put(0.2,0.3){\footnotesize $7$}
\end{picture}, 
\begin{picture}(1,5)
\put(0,0){\line(0,1){5}} 
\put(0.7,0){\line(0,1){5}} 
\put(0,0){\line(1,0){0.7}} 
\put(0,1){\line(1,0){0.7}} 
\put(0,2){\line(1,0){0.7}} 
\put(0,3){\line(1,0){0.7}} 
\put(0,4){\line(1,0){0.7}} 
\put(0,5){\line(1,0){0.7}} 
\put(0.2,4.3){\footnotesize $1$}
\put(0.2,3.3){\footnotesize $2$}
\put(0.2,2.3){\footnotesize $3$}
\put(0.2,1.3){\footnotesize $6$}
\put(0.2,0.3){\footnotesize $7$}
\end{picture}, 
\begin{picture}(1,5)
\put(0,0){\line(0,1){5}} 
\put(0.7,0){\line(0,1){5}} 
\put(0,0){\line(1,0){0.7}} 
\put(0,1){\line(1,0){0.7}} 
\put(0,2){\line(1,0){0.7}} 
\put(0,3){\line(1,0){0.7}} 
\put(0,4){\line(1,0){0.7}} 
\put(0,5){\line(1,0){0.7}} 
\put(0.2,4.3){\footnotesize $1$}
\put(0.2,3.3){\footnotesize $2$}
\put(0.2,2.3){\footnotesize $4$}
\put(0.2,1.3){\footnotesize $5$}
\put(0.2,0.3){\footnotesize $7$}
\end{picture}, 
\begin{picture}(1,5)
\put(0,0){\line(0,1){5}} 
\put(0.7,0){\line(0,1){5}} 
\put(0,0){\line(1,0){0.7}} 
\put(0,1){\line(1,0){0.7}} 
\put(0,2){\line(1,0){0.7}} 
\put(0,3){\line(1,0){0.7}} 
\put(0,4){\line(1,0){0.7}} 
\put(0,5){\line(1,0){0.7}} 
\put(0.2,4.3){\footnotesize $1$}
\put(0.2,3.3){\footnotesize $2$}
\put(0.2,2.3){\footnotesize $4$}
\put(0.2,1.3){\footnotesize $6$}
\put(0.2,0.3){\footnotesize $7$}
\end{picture}, 
\begin{picture}(1,5)
\put(0,0){\line(0,1){5}} 
\put(0.7,0){\line(0,1){5}} 
\put(0,0){\line(1,0){0.7}} 
\put(0,1){\line(1,0){0.7}} 
\put(0,2){\line(1,0){0.7}} 
\put(0,3){\line(1,0){0.7}} 
\put(0,4){\line(1,0){0.7}} 
\put(0,5){\line(1,0){0.7}} 
\put(0.2,4.3){\footnotesize $1$}
\put(0.2,3.3){\footnotesize $2$}
\put(0.2,2.3){\footnotesize $5$}
\put(0.2,1.3){\footnotesize $6$}
\put(0.2,0.3){\footnotesize $7$}
\end{picture}, 
\begin{picture}(1,5)
\put(0,0){\line(0,1){5}} 
\put(0.7,0){\line(0,1){5}} 
\put(0,0){\line(1,0){0.7}} 
\put(0,1){\line(1,0){0.7}} 
\put(0,2){\line(1,0){0.7}} 
\put(0,3){\line(1,0){0.7}} 
\put(0,4){\line(1,0){0.7}} 
\put(0,5){\line(1,0){0.7}} 
\put(0.2,4.3){\footnotesize $1$}
\put(0.2,3.3){\footnotesize $3$}
\put(0.2,2.3){\footnotesize $4$}
\put(0.2,1.3){\footnotesize $5$}
\put(0.2,0.3){\footnotesize $7$}
\end{picture}, 
\begin{picture}(1,5)
\put(0,0){\line(0,1){5}} 
\put(0.7,0){\line(0,1){5}} 
\put(0,0){\line(1,0){0.7}} 
\put(0,1){\line(1,0){0.7}} 
\put(0,2){\line(1,0){0.7}} 
\put(0,3){\line(1,0){0.7}} 
\put(0,4){\line(1,0){0.7}} 
\put(0,5){\line(1,0){0.7}} 
\put(0.2,4.3){\footnotesize $1$}
\put(0.2,3.3){\footnotesize $3$}
\put(0.2,2.3){\footnotesize $4$}
\put(0.2,1.3){\footnotesize $6$}
\put(0.2,0.3){\footnotesize $7$}
\end{picture}, 
\begin{picture}(1,5)
\put(0,0){\line(0,1){5}} 
\put(0.7,0){\line(0,1){5}} 
\put(0,0){\line(1,0){0.7}} 
\put(0,1){\line(1,0){0.7}} 
\put(0,2){\line(1,0){0.7}} 
\put(0,3){\line(1,0){0.7}} 
\put(0,4){\line(1,0){0.7}} 
\put(0,5){\line(1,0){0.7}} 
\put(0.2,4.3){\footnotesize $1$}
\put(0.2,3.3){\footnotesize $3$}
\put(0.2,2.3){\footnotesize $5$}
\put(0.2,1.3){\footnotesize $6$}
\put(0.2,0.3){\footnotesize $7$}
\end{picture}, 
\begin{picture}(1,5)
\put(0,0){\line(0,1){5}} 
\put(0.7,0){\line(0,1){5}} 
\put(0,0){\line(1,0){0.7}} 
\put(0,1){\line(1,0){0.7}} 
\put(0,2){\line(1,0){0.7}} 
\put(0,3){\line(1,0){0.7}} 
\put(0,4){\line(1,0){0.7}} 
\put(0,5){\line(1,0){0.7}} 
\put(0.2,4.3){\footnotesize $1$}
\put(0.2,3.3){\footnotesize $4$}
\put(0.2,2.3){\footnotesize $5$}
\put(0.2,1.3){\footnotesize $6$}
\put(0.2,0.3){\footnotesize $7$}
\end{picture}, 
\begin{picture}(1,5)
\put(0,0){\line(0,1){5}} 
\put(0.7,0){\line(0,1){5}} 
\put(0,0){\line(1,0){0.7}} 
\put(0,1){\line(1,0){0.7}} 
\put(0,2){\line(1,0){0.7}} 
\put(0,3){\line(1,0){0.7}} 
\put(0,4){\line(1,0){0.7}} 
\put(0,5){\line(1,0){0.7}} 
\put(0.2,4.3){\footnotesize $2$}
\put(0.2,3.3){\footnotesize $3$}
\put(0.2,2.3){\footnotesize $4$}
\put(0.2,1.3){\footnotesize $5$}
\put(0.2,0.3){\footnotesize $7$}
\end{picture}, 
\begin{picture}(1,5)
\put(0,0){\line(0,1){5}} 
\put(0.7,0){\line(0,1){5}} 
\put(0,0){\line(1,0){0.7}} 
\put(0,1){\line(1,0){0.7}} 
\put(0,2){\line(1,0){0.7}} 
\put(0,3){\line(1,0){0.7}} 
\put(0,4){\line(1,0){0.7}} 
\put(0,5){\line(1,0){0.7}} 
\put(0.2,4.3){\footnotesize $2$}
\put(0.2,3.3){\footnotesize $3$}
\put(0.2,2.3){\footnotesize $4$}
\put(0.2,1.3){\footnotesize $6$}
\put(0.2,0.3){\footnotesize $7$}
\end{picture}, 
\begin{picture}(1,5)
\put(0,0){\line(0,1){5}} 
\put(0.7,0){\line(0,1){5}} 
\put(0,0){\line(1,0){0.7}} 
\put(0,1){\line(1,0){0.7}} 
\put(0,2){\line(1,0){0.7}} 
\put(0,3){\line(1,0){0.7}} 
\put(0,4){\line(1,0){0.7}} 
\put(0,5){\line(1,0){0.7}} 
\put(0.2,4.3){\footnotesize $2$}
\put(0.2,3.3){\footnotesize $3$}
\put(0.2,2.3){\footnotesize $5$}
\put(0.2,1.3){\footnotesize $6$}
\put(0.2,0.3){\footnotesize $7$}
\end{picture}, 
\begin{picture}(1,5)
\put(0,0){\line(0,1){5}} 
\put(0.7,0){\line(0,1){5}} 
\put(0,0){\line(1,0){0.7}} 
\put(0,1){\line(1,0){0.7}} 
\put(0,2){\line(1,0){0.7}} 
\put(0,3){\line(1,0){0.7}} 
\put(0,4){\line(1,0){0.7}} 
\put(0,5){\line(1,0){0.7}} 
\put(0.2,4.3){\footnotesize $2$}
\put(0.2,3.3){\footnotesize $4$}
\put(0.2,2.3){\footnotesize $5$}
\put(0.2,1.3){\footnotesize $6$}
\put(0.2,0.3){\footnotesize $7$}
\end{picture}, and 
\begin{picture}(1,5)
\put(0,0){\line(0,1){5}} 
\put(0.7,0){\line(0,1){5}} 
\put(0,0){\line(1,0){0.7}} 
\put(0,1){\line(1,0){0.7}} 
\put(0,2){\line(1,0){0.7}} 
\put(0,3){\line(1,0){0.7}} 
\put(0,4){\line(1,0){0.7}} 
\put(0,5){\line(1,0){0.7}} 
\put(0.2,4.3){\footnotesize $3$}
\put(0.2,3.3){\footnotesize $4$}
\put(0.2,2.3){\footnotesize $5$}
\put(0.2,1.3){\footnotesize $6$}
\put(0.2,0.3){\footnotesize $7$}
\end{picture},
\end{center}
a total of 14 columnar tableaux.  
This count agrees with the $(n,k)=(7,5)$ entry of the example array $\mbox{\Large (}\mya_{4}(n,k)\mbox{\Large )}$ depicted in \S 2 above. 
At this point, a correspondence with the fourteen $4$-admissible $(7,5)$-NE-paths presented above should be clear. 
Now let 
\begin{equation}
\myc'_{t}(n,k) := \mbox{\Large $|$}\{\mbox{\small 
$t$-admissible $(n,k)$-columnar tableaux}\}\mbox{\Large $|$}.
\end{equation}
The proof of our next lemma is obtained by producing an explicit bijection between the $t$-admissible $(n,k)$-columnar tableaux and the $t$-admissible $(n,k)$-NE-paths. 
This bijection merely formalizes what we have observed in our example correspondence between the $4$-admissible $(7,5)$-NE-paths and the $4$-admissible $(7,5)$-columnar tableaux.

\noindent 
{\bf \LemmaOnePrime}\ \ {\sl For all integers $n$ and $k$, we have $\myc_{t}(n,k) = \myc'_{t}(n,k)$.}

{\em Proof.} Let $\mathcal{P}_{t}(n,k)$ be the set of $t$-admissible $(n,k)$-NE-paths, and let $\mathcal{T}_{t}(n,k)$ be the set of $t$-admissible $(n,k)$-columnar tableaux. 
In this proof, we identify an $(n,k)$-NE-path $\selt$ with the sequence $\selt = \left(\rule[-2.5mm]{0mm}{6mm}(x_{i}(\selt),y_{i}(\selt))\right)_{i=1}^{k}$ consisting of the $k$ successive endpoints of the horizontal, or easterly, steps of the path. 
For example, for the first $4$-admissible $(7,5)$-NE-path depicted above, the sequence of horizontal endpoints is $\left(\rule[-2.5mm]{0mm}{6mm}(1,0),(2,0),(3,0),(4,1),(5,2)\right)$. 

Given a $t$-admissible $(n,k)$-NE-path $\selt = \left(\rule[-2.5mm]{0mm}{6mm}(x_{i}(\selt),y_{i}(\selt))\right)_{i=1}^{k}$, set $\phi(\selt) := \left(\rule[-2.5mm]{0mm}{6mm}x_{i}(\selt)+y_{i}(\selt)\right)_{i=1}^{k}$. 
Let $T = (T_{1},\ldots,T_{k}) = \phi(\selt)$. 
Since $x_{i}(\selt)=i$, then $i \leq x_{i}(\selt)+y_{i}(\selt) = T_{i}$. 
In particular, $1 \leq T_{1}$. 
Also, $y_{i}(\selt) \leq y_{i+1}(\selt)$, then $x_{i}(\selt)+y_{i}(\selt) = T_{i} < T_{i+1} = x_{i+1}(\selt)+y_{i+1}(\selt)$ when $i \in \{1,2,\ldots,k-1\}$.  
Since $x_{k}(\selt)=k$ and $y_{k}(\selt) \leq n-k$, then $T_{k} \leq n$.
So $T$ is an $(n,k)$-columnar tableau. 
Now suppose $t-1<k$, and let $j \in \{1,2,\ldots,k-t+1\}$. 
Since $\selt$ is $t$-admissible, then we have $j = (t-1+j)-t+1 = x_{t-1+j}(\selt)-t+1 \leq y_{t-1+j}(\selt)$. 
So, $t-1+2j = (t-1+j) +j = x_{t-1+j}(\selt) + j \leq x_{t-1+j}(\selt)+y_{t-1+j}(\selt) = T_{t-1+j}$. 
Thus $T$ is a $t$-admissible $(n,k)$-columnar tableau. 
We can therefore regard $\phi: \mathcal{P}_{t}(n,k) \longrightarrow \mathcal{T}_{t}(n,k)$ as a well-defined function. 

Now suppose that $T = (T_{1},\ldots,T_{k})$ is a $t$-admissible $(n,k)$-columnar tableau. 
Declare that $\psi(T) := \left(\rule[-2.5mm]{0mm}{6mm}(i,T_{i}-i)\right)_{i=1}^{k}$. 
Set $\selt = \left(\rule[-2.5mm]{0mm}{6mm}(x_{i},y_{i})\right)_{i=1}^{k} := \psi(T)$. 
To prove that $\selt$ is an $(n,k)$-NE-path, it suffices to check that $0 \leq y_{1} \leq y_{2} \leq \cdots \leq y_{k} \leq n-k$. 
Since $1 \leq T_{1}$, then $0 \leq T_{1}-1 = y_{1}$. 
When $i \in \{1,2,\ldots,k-1\}$, then $T_{i} < T_{i+1}$ so $y_{i} = T_{i}-i \leq T_{i+1}-(i+1) = y_{i+1}$. 
Also, $y_{k} = T_{k}-k \leq n-k$ since $T_{k} \leq n$. 
Next, we check that $\selt$ is $t$-admissible by showing that each $y_{i} \geq x_{i}-t+1$, assuming $t-1 < k$. 
Suppose for the moment that $i>t-1$, so that $i=t-1+j$ with $j \in \{1,2,\ldots,k-t+1\}$. 
Then $y_{i} = T_{i}-i = T_{t-1+j}-(t-1+j) \geq t-1+2j-t+1-j = j = i-t+1 = x_{i}-t+1$. 
On the other hand, if $i \leq t-1$, then $x_{i}-t+1 = i-t+1 \leq 0 \leq y_{i}$. 
This reasoning shows that $\selt$ is $t$-admissible. 

We can therefore regard $\psi: \mathcal{T}_{t}(n,k) \longrightarrow \mathcal{P}_{t}(n,k)$ as a well-defined function.
Clearly $\phi$ and $\psi$ are inverses, so $\mathcal{P}_{t}(n,k)$ and $\mathcal{T}_{t}(n,k)$ are equinumerous, which is what we needed to show.\hfill\QED

Finally, we consider another set of numbers $\myd_{t}(n,k)$ which arise as 
coefficients of certain polynomials or, from another viewpoint, as 
entries of representing matrices for certain linear transformations on 
polynomial vector spaces. 
It appears this can be viewed within the context of Rota's finite operator calculus (see \cite{Niederhausen}), but we use more direct and elementary reasoning here. 
To set things up, let $\{x_{j}\}_{j \geq 0}$ 
be the basis for the polynomial vector space $\mathbb{R}[x]$ 
(polynomials in the indeterminate $x$ and with real coefficients) given 
by $x_{j} := x^{j}/j!$.  Let $S$ be the subspace of $\mathbb{R}[x]$ spanned by 
$\{x_{j}\}_{j \geq 1}$. The linear transformation $D: S 
\longrightarrow \mathbb{R}[x]$ will be the differential operator $D(y) := 
y' + y''$.  In fact, $D$ is a vector space isomorphism.  For any 
positive integer $N$, set $D^{-N} := (D^{-1})^{N}$. 

It is easy to see by induction on $N$ that for all positive integers 
$N$, $D^{-N}(x_{t-1})$ is in 
$\mbox{span}_{\mathbb{R}}\{x_{j}\}_{j=1}^{t-1+N}$.  (For the basis step 
of the induction argument, check that $\displaystyle D^{-1}(x_{t-1}) = 
\sum_{i=0}^{t-1}(-1)^{i}x_{t-i}$.) This means that  for any positive 
integer $N$ we can write  
\begin{equation}
D^{-N}(x_{t-1}) = \sum_{i=0}^{t-1+N-1}(-1)^{i} \myd_{t}(i+N-1,i) x_{t-1+N-i}
\end{equation}
for some real numbers $\myd_{t}(i+N-1,i)$. Declare $\myd_{t}(n,k)$ to be 
zero for any integer pair $(n,k)$ such that for all $N > 0$ and $0 
\leq i \leq t-1+N-1$, $(n,k) \not= (i+N-1,i)$, i.e.\ $\myd_{t}(n,k) = 0$ if $(n,k)$ does not 
index any term  
appearing in any of the sums shown in equation (4) above when $N \geq 1$. 
Thus, $\myd_{t}$ is a function defined on all of $\mathbb{Z} \times \mathbb{Z}$. 

\noindent 
{\bf \LemmaTwo}\ \ {\sl (i) We have $\myd_{t}(0,0) = 1$.  For all integers 
$n$ and $k$, we have (ii) $\myd_{t}(n,k) = 0$ if $n<0$, 
$k<0$, or $k > \min\left\{\lfloor \frac{n-1+t}{2} \rfloor , 
n\right\}$, and otherwise 
(iii) $\myd_{t}(n,k) = \myd_{t}(n-1,k-1) + \myd_{t}(n-1,k)$ as long as $n>0$.} 

{\em Proof.} For ({\sl i}), consider the $i=0$ term in the expression 
for $D^{-1}(x_{t-1})$ given in the paragraph preceding the lemma 
statement.  Then $1 = \myd_{t}(i+N-1,i) = \myd_{t}(0,0)$.  For ({\sl ii}), we observe that 
an integer pair 
$(n,k)$ is indeed a pair $(i+N-1,i)$ corresponding to a term 
in the sum shown in 
equation (4) above if and only if $k=i$ and $n = k+N-1$ for some $N>0$ and $0 \leq i \leq t-1+N-1$.  Now simply check inequalities to see that $k \geq 0$, $n \geq 0$, $k \leq n$, and $k \leq \lfloor \frac{n-1+t}{2} \rfloor$. 


For ({\sl iii}), we apply $D$ to each side of 
equation (4).  First, 
\begin{equation}
D(D^{-N}(x_{t-1})) = D^{-(N-1)}(x_{t-1}) = 
\sum_{i=0}^{t-1+N-2}(-1)^{i} \myd_{t}(i+N-2,i) x_{t-1+N-1-i}.
\end{equation} 
On the other hand, 
\begin{eqnarray*}
D(D^{-N}(x_{t-1})) & = & D\bigg(\sum_{i=0}^{t-1+N-1}(-1)^{i} \myd_{t}(i+N-1,i) 
x_{u+N-i}\bigg)\\ 
& = & \sum_{i=0}^{t-1+N-1}(-1)^{i} \myd_{t}(i+N-1,i) 
(x_{t-1+N-1-i}+x_{t-1+N-2-i})\\ 
& = & \sum_{i=0}^{t-1+N-1}(-1)^{i} [\myd_{t}(i+N-1,i)-\myd_{t}(i+N-2,i-1)] 
x_{t-1+N-1-i}, 
\end{eqnarray*}
where the latter is obtained by expanding and reindexing. Then 
by equating coefficients in the latter expression with coefficients 
for the expression obtained in equation (5), 
we see that $\myd_{t}(t-1+2N-2,t-1+N-1) - 
\myd_{t}(t-1+2N-3,t-1+2N-2) = 0$ and that 
for all $0 \leq i \leq t-1+N-2$ we have 
$\myd_{t}(i+N-1,i)-\myd_{t}(i+N-2,i-1) = \myd_{t}(i+N-2,i)$. The latter formula 
actually becomes the former  
when $i = t-1+N-1$, as $\myd_{t}(t-1+2N-3,t-1+N-1)$ evaluates to zero. 
Therefore, 
\begin{equation}
\myd_{t}(i+N-1,i) = \myd_{t}(i+N-2,i-1) + \myd_{t}(i+N-2,i) 
\end{equation} for 
all $N>0$ and $0 \leq i \leq t-1+N-1$.  
Now if $(n,k)$ is an integer pair with $n \geq 0$ and $0 \leq k 
\leq \min\left\{\lfloor \frac{n-1+t}{2} \rfloor , 
n\right\}$, then set $i=k$ and $N=n+1-k$.  As in the previous 
paragraph we can see that $(i+N-1,i)$ corresponds to a term from 
equation (4).  Then from equation (6), we get $\myd_{t}(n,k) = 
\myd_{t}(n-1,k-1) + \myd_{t}(n-1,k)$, as desired.\hfill\QED 

The preceding lemmas furnish the key steps in the proof of the main result of this section, whose value is not so much the novelty of the results (which are routine to enumeration experts) but rather the pleasantness and illustrative utility of the results taken together as an assemblage.  

\noindent 
{\bf \MainTheorem}\ \ {\sl For all integers $n$ and $k$, we have}
\[\mya_{t}(n,k) = \myb_{t}(n,k) = \myc_{t}(n,k) = \myc'_{t}(n,k) = \myd_{t}(n,k).\]

{\em Proof.} The $\myb_{t}(n,k)$'s are easily seen to satisfy 
the defining conditions ({\sl i}), ({\sl ii}), and ({\sl iii}) which uniquely 
determine the $\mya_{t}(n,k)$'s, so $\myb_{t}(n,k) = \mya_{t}(n,k)$.  
The $\myd_{t}(n,k)$'s satisfy these 
same conditions by \LemmaTwo, hence $\myd_{t}(n,k) = \mya_{t}(n,k)$.  
And by \BothLemmaOnes\ we have $\myb_{t}(n,k) = \myc_{t}(n,k) = \myc'_{t}(n,k)$.\hfill\QED

\vspace*{0.1in}
{\bf \S 4\ \  Counting odd numbers in truncations of Pascal's triangles.} 
It is a well-known phenomenon that the number of odds on any given row of Pascal's triangle is a power of two. 
A classical proof of this fact utilizes Lucas's Theorem and is recapitulated in \LucasCorollary\ below, cf.\ \S 8.4 of \cite{BQ}.
Before we proceed, we fix some notation. 
For a prime $p$ and a nonnegative integer $m$, let $l_{p}(m)$ be $0$ when $m=0$ and $\lfloor \log_{p}(m) \rfloor$ otherwise. 
For $i \in \{0,1,\ldots,l_{p}(m)\}$, let $m^{(p)}_{i}$ denote the $i^{\mbox{\tiny th}}$ digit of the base $p$ representation of $m$, and let $\mathcal{D}_{p}(m)$ be the subset of $\{0,1,\ldots,l_{p}(m)\}$ for which $i \in \mathcal{D}_{p}(m)$ if and only if $m^{(p)}_{i} \ne 0$. 
Further, let $d_{p}(n) := |\mathcal{D}_{p}(m)|$. 
Our interest is mainly in the case that $p=2$, but we state Lucas's Theorem in its full generality in order to encourage the reader (and the authors) to keep in mind the possibility of extending some of the ideas of this section to other primes. 

\noindent 
{\bf \LucasTheorem\ (Lucas's Theorem)}\ \ {\sl Let $p$ be any prime, and fix any nonnegative integers $n$ and $k$. Let} $l := \mbox{max}(l_{p}(n),l_{p}(k))$.  {\sl Then,}
\[{n \choose k} \equiv \prod_{i=0}^{l}{n^{(p)}_{i} \choose k^{(p)}_{i}}\ \mbox{$($mod } p).\]

From here on, we suppress the ``$p$'' superscripts and subscripts from the notation introduced above and fix $p=2$ as our prime. 
So, for example, ``$d(n)$'' means $d_{2}(n)$, ``$l(n)$'' means $l_{2}(n)$, ``$\mathcal{D}(n)$'' means $\mathcal{D}_{2}(n)$, etc. 

\noindent 
{\bf \LucasCorollary}\ \ {\sl Let $n$ be a nonnegative integer. Then the number of odds on the} $n^{\mbox{\tiny th}}$ {\sl row of Pascal's triangle is $2^{d(n)}$.} 

{\em Proof.} Suppose $0 \leq k \leq n$. 
Then by Lucas' Theorem, ${n \choose k}$ is odd if and only if $n_{i} = 1$ whenever $k_{i}=1$. 
So, ${n \choose k}$ is odd if and only if $\mathcal{D}(k) \subseteq \mathcal{D}(n)$. 
Of course, there are $2^{|\mathcal{D}(n)|} = 2^{d(n)}$ choices for such subsets.\hfill\QED

Now we turn our attention to truncations of Pascal's triangle. 
\LucasCorollary\ and the reasoning exhibited in its proof will be used in several of the lemmas that follow.  
These lemmas support the proof of the following theorem, which is the main result of this paper. 

\noindent 
{\bf \OddCountingTheorem}\ \ {\sl The number of odds on each row of the Pascal triangle truncation} $\mbox{\Large (}\mya_{t}(n,k)\mbox{\Large )}$ {\sl is a power of two if and only if $t$ is a power of two. 
In this case, when $n$ is a nonnegative integer, the number of odds on row $n$ of the array is precisely $\frac{1}{2} \cdot 2^{d(n+t)}$.}

The proof of \OddCountingTheorem\ is at the end of this section and will be easily deduced from the lemmas we establish next. 

\noindent 
{\bf \CountingLemmaConverse}\ \ {\sl Suppose $t$ is not a power of two.  Then the number of odds on row $t$ of the truncated Pascal's triangle} $\mbox{\Large (}\mya_{t}(n,k)\mbox{\Large )}$ {\sl is an odd number greater than one and therefore not a power of two.} 

{\em Proof.}  Let $n=t$. 
If $0 \leq k \leq n-1$, then $\mya_{t}(n,k) = {n \choose k} - {n \choose k-t} = {n \choose k} - {n \choose k-n} = {n \choose k}$. 
And if $k=n$, then $\mya_{t}(n,k) = {n \choose n} - {n \choose 0} = 0$. 
So the $n^{\mbox{\tiny th}}$ row of the given truncated Pascal array is the same as the $n^{\mbox{\tiny th}}$ row of Pascal's triangle with the sole exception of the $n^{\mbox{\tiny th}}$ entry, which is a 1 in Pascal's triangle and a 0 in the truncated Pascal array. 
So, the number of odds on row $n$ of the truncated Pascal array is, by \LucasCorollary, $2^{d(n)}-1$. 
This odd number is a power of two if and only if $d(n)=0$ if and only if $n=0$ if and only if $t=0$.  
Since $t$ is a positive integer, we conclude that the number of odds on row $n$ is an odd number greater than one.\hfill\QED

The simple observations of \CountingLemmaS\ are needed for our proof of \CountingLemmaParity. 

\noindent 
{\bf \CountingLemmaCentral}\ \ {\sl Let $m$ be a nonnegative integer. 
The binomial coefficient ${2(m+1) \choose m+1}$ is even. 
The binomial coefficient ${2m+1 \choose m}$ is odd if and only if there is a positive integer $q$ such that $2m+1=2^{q}-1$.} 

{\em Proof.} Since ${2(m+1) \choose m+1} = {2m+1 \choose m+1} + {2m+1 \choose m}$, which is even since ${2m+1 \choose m+1} = {2m+1 \choose m}$. 
Now assume $2m+1 =  2^{q}-1$ for a positive integer $q$. 
If $q=1$, then ${2m+1 \choose m} = {1 \choose 0} = 1$, which is odd. 
Next assume $q > 1$. 
Now, $l(2m+1) = q-1$, $l(m)=q-2$, $\mathcal {D}(2m+1) = \{0,1,\ldots,q-1\}$, and $\mathcal {D}(m) = \{0,1,\ldots,q-2\}$. 
So by Lucas's Theorem, ${2m+1 \choose m} \equiv {1 \choose 0}{1 \choose 1}{1 \choose 1} \cdots {1 \choose 1} \mbox{ $($mod  $2)$}$, hence ${2m+1 \choose m}$ is odd. 
Finally, assume ${2m+1 \choose m}$ is odd. 
If $m=0$, then $2m+1=1=2^{1}-1$. 
So now assume $m > 0$. 
Set $r-2 := l(m)$ (hence $r > 0$) and write $m = m_{r-2}2^{r-2} + m_{r-3}2^{r-3} + \cdots + m_{1}2^{1} + m_{0}2^{0}$, where of course $m_{r-2} = 1$. 
So, $2m+1 = m_{r-2}2^{r-1} + m_{r-3}2^{r-2} + \cdots + m_{0}2^{1} + 1 \cdot 2^{0}$. 
Since ${2m+1 \choose m}$ is odd, Lucas's Theorem requires that ${m_{r-3} \choose m_{r-2}} = {m_{r-4} \choose m_{r-3}} = \cdots = {m_{0} \choose m_{1}} = 1$. 
Based on these equalities, we observe that $m_{r-2}=1$ forces $m_{r-3} = 1$, which in turn forces $m_{r-4}=1$, etc. 
We conclude that $m_{r-2} = m_{r-3} = m_{r-4} = \cdots = m_{1} = m_{0} = 1$. 
Then $m=2^{r-1}-1$, so $2m+1=2^{r}-1$.\hfill\QED

\noindent 
{\bf \CountingLemmaOddRow}\ \ {\sl Let $n$ be a nonnegative integer.  All entries on the} $n^{\mbox{\tiny th}}$ {\sl row of Pascal's triangle are odd if and only if there is a nonnegative integer $q$ such that $n=2^{q}-1$.} 

{\em Proof.} Suppose all entries on the $n^{\mbox{\tiny th}}$ row are odd. 
If $n=0$, then $n = 2^{0}-1$. 
If $n>0$, then by \CountingLemmaCentral, there is a positive integer $q$ with $n=2^{q}-1$. 
Conversely, suppose $n=2^{q}-1$ for a nonnegative integer $q$. 
If $n=0=2^{0}-1$, then all entries on this row are odd, since the only entry on this row is ${0 \choose 0} = 1$. 
Now say $q$ is positive, so $n = n_{q-1}2^{q-1} + n_{q-2}2^{q-2} + \cdots + n_{1}2^{1} + n_{0}2^{0}$ with $n_{q-1} = n_{q-2} = \cdots = n_{1} = n_{0} = 1$. 
Let $0 \leq k \leq n$, and write $k = k_{q-1}2^{q-1} + k_{q-2}2^{q-2} + \cdots + k_{1}2^{1} + k_{0}2^{0}$. 
Then by Lucas's Theorem, $\displaystyle {n \choose k} \equiv {1 \choose k_{q-1}}{1 \choose k_{q-2}} \cdots {1 \choose k_{0}}\ \mbox{$($mod $2)$} \equiv 1 \ \mbox{$($mod $2)$}$. 
So all entries on the $n^{\mbox{\tiny th}}$ row are odd.\hfill\QED

The following binomial coefficient identity is a version of Vandermonde's Identity, cf.\ Identity 132 of \cite{BQ}. 

\noindent 
{\bf \CountingLemmaIdentity\ (Vandermonde's Identity)}\ \ {\sl Let $m$, $l$, and $j$ be nonnegative integers. Then}
\[\sum_{i=0}^{l}{l \choose i}{m-l \choose j-i} = {m \choose j}.\]

\noindent 
{\bf \CountingLemmaParity}\ \ {\sl Suppose $t=2^q$ for some nonnegative integer $q$.  Then the quantities} $\mya_{t}(n,k)$ {\sl and ${n+t \choose k}$ have the same parity.}

{\em Proof.} By \CountingLemmaIdentity, ${n+t \choose k} = \sum_{i=0}^{t}{t \choose i}{n+t-t \choose k-i}$. 
Then ${n+t \choose k} \equiv \sum_{i=0}^{t}{t \choose i}{n \choose k-i}\ \mbox{$($mod $2)$}$. 
Since $t = 2^{q}$ for a nonnegative integer $q$, then by \CountingLemmaOddRow, all entries on row $t-1$ are odd. 
Then all entries on row $t$ except the first and last are even. 
So, $\sum_{i=0}^{t}{t \choose i}{n \choose k-i} \equiv {n \choose k} + {n \choose k-t}\ \mbox{$($mod $2)$}$. 
And, ${n \choose k} + {n \choose k-t} \equiv {n \choose k} - {n \choose k-t}\ \mbox{$($mod $2)$}$. 
Since $\mya_{t}(n,k) = {n \choose k} - {n \choose k-t}$ by \MainTheorem, we conclude that ${n+t \choose k} \equiv \mya_{t}(n,k)\ \mbox{$($mod $2)$}$.\hfill\QED

We are now ready to prove \OddCountingTheorem. 

{\em Proof of \OddCountingTheorem.} \CountingLemmaConverse\ shows that if the number of odds on each row of the Pascal triangle truncation $\mbox{\Large (}\mya_{t}(n,k)\mbox{\Large )}$ is a power of two, then $t$ must be a power of two. 
Conversely, let us now suppose that $t = 2^{q}$ for some nonnegative integer $q$. 
We aim to demonstrate the following claim: When $n$ is a nonnegative integer, the number of odds on row $n$ of the array is precisely $\frac{1}{2} \cdot 2^{d(n+t)}$. 

We begin by assuming $n$ is odd. 
Write $n=2m+1$ for a nonnegative integer $m$. 
The last nonzero entry on row $n$ occurs at position $k = \mbox{min}\left\{\left\lfloor\frac{t-1+n}{2}\right\rfloor,n\right\} = \mbox{min}\left\{\left\lfloor\frac{2^{q}+2m}{2}\right\rfloor,2m+1\right\}$. 
If $q=0$, then $\left\lfloor\frac{2^{q}+2m}{2}\right\rfloor = m$, so $k=m$. 
Of course, we now have $t=2^{0}=1$. 
By \CountingLemmaParity, the parity of entry $\mya_{1}(n,j)$ of the $n^{\mbox{\tiny th}}$ row of our array (where $0 \leq j \leq k=m$) is the same as the parity of the binomial coefficient ${n+1 \choose j} = {2m+2 \choose j}$. 
Since ${2m+2 \choose m+1}$ is even by \CountingLemmaCentral, then the number of odds on the $(n+1)^{\mbox{\tiny st}}$ row of Pascal's triangle is twice the number of odds amongst the entries entry ${n+1 \choose j}$ for $0 \leq j \leq k=m$. 
Now, the number of odds on the $(n+1)^{\mbox{\tiny st}}$ row of Pascal's triangle is $2^{d(n+1)}$ by \LucasCorollary. 
Therefore, the number of odds on the $n^{\mbox{\tiny th}}$ row of our truncated Pascal array is $\frac{1}{2} \cdot 2^{d(n+1)}$, confirming our desired claim when $q=0$. 

Continuing with the assumption that $n$ is odd, now assume that $q>0$. 
Then $\left\lfloor\frac{2^{q}+2m}{2}\right\rfloor = 2^{q-1}+m$. 
If $2^{q-1}+m \geq 2m+1$ (and hence $2^{q} > 2m+1$), then $k=2m+1$, so the entries of the $n^{\mbox{\tiny th}}$ row of our array coincide with the entries of the $n^{\mbox{\tiny th}}$ row of Pascal's triangle. 
This shared number of odds is therefore $2^{d(n)}$, by \LucasCorollary. 
But $d(n+t) = d(n+2^{q}) =d(n)+1$ since $2^{q} > 2m+1 = n$. 
Then the number of odds on the $n^{\mbox{\tiny th}}$ row of our array is $2^{d(n)} = \frac{1}{2} \cdot 2^{d(n)+1} = \frac{1}{2} \cdot 2^{d(n+t)}$, again confirming our claim. 
So now consider the case that  $2^{q-1}+m < 2m+1$. 
The $n^{\mbox{\tiny th}}$ row entry $\mya_{t}(n,j)$ of our truncated Pascal array (where $0 \leq j \leq k$) has the same parity as the entry ${n+t \choose j}$ of Pascal's triangle, by  \CountingLemmaParity. 
Now, $k=2^{q-1}+m$ while $n+t = 2^{q}+2m+2 = 2(2^{q-1}+m+1)$. 
Since ${n+t \choose 2^{q-1}+m+1} = {2(2^{q-1}+m+1) \choose 2^{q-1}+m+1}$ is even by \CountingLemmaCentral, then the number of odds on the $(n+t)^{\mbox{\tiny th}}$ row of Pascal's triangle is twice the number of odds amongst the entries entry ${n+t \choose j}$ for $0 \leq j \leq k = 2^{q-1}+m$. 
So, the number of odds on the $n^{\mbox{\tiny th}}$ row of our truncated Pascal array is $\frac{1}{2} \cdot 2^{d(n+t)}$, completing the confirmation of our claim when $n$ is odd. 

Next, assume $n$ is even, and write $n=2m$ for some nonnegative integer $m$. 
As before, the last nonzero entry on row $n$ occurs at position $k = \mbox{min}\left\{\left\lfloor\frac{t-1+n}{2}\right\rfloor,n\right\} = \mbox{min}\left\{\left\lfloor\frac{2^{q}-1+2m}{2}\right\rfloor,2m\right\}$. 
Say $q=0$, so $t=1$. 
Then $k = \mbox{min}\left\{\left\lfloor\frac{2m}{2}\right\rfloor,2m\right\} = m$. 
By \CountingLemmaParity, the parity of entry $\mya_{1}(n,j)$ of the $n^{\mbox{\tiny th}}$ row of our array (where $0 \leq j \leq k=m$) is the same as the parity of the binomial coefficient ${n+1 \choose j} = {2m+1 \choose j}$. 
Since the entries ${n+1 \choose 0}, {n+1 \choose 1}, \ldots, {n+1 \choose m}$ comprise exactly half of the entries of said row of Pascal's triangle, then there are $\frac{1}{2} \cdot 2^{d(n+1)}$ odds amongst these entries. 
So there are  $\frac{1}{2} \cdot 2^{d(n+1)}$ odd entries on the $n^{\mbox{\tiny th}}$ row of our truncated Pascal array, completing the confirmation of our claim when $q=0$.  

Keeping the hypothesis that $n$ is even, now assume that $q>0$. 
Then $\left\lfloor\frac{2^{q}-1+2m}{2}\right\rfloor = 2^{q-1}+m-1$. 
If $2^{q-1}+m-1 \geq 2m$ (and hence $2^{q} > 2m$), then $k=2m$, so the entries of the $n^{\mbox{\tiny th}}$ row of our array coincide with the entries of the $n^{\mbox{\tiny th}}$ row of Pascal's triangle. 
This shared number of odds is therefore $2^{d(n)}$, by \LucasCorollary. 
But $d(n+t) = d(n+2^{q}) =d(n)+1$ since $2^{q} > 2m = n$. 
Then the number of odds on the $n^{\mbox{\tiny th}}$ row of our array is $2^{d(n)} = \frac{1}{2} \cdot 2^{d(n)+1} = \frac{1}{2} \cdot 2^{d(n+t)}$, again confirming our claim. 
So now consider the case that  $2^{q-1}+m-1 < 2m$. 
The $n^{\mbox{\tiny th}}$ row entry $\mya_{t}(n,j)$ of our truncated Pascal array (where $0 \leq j \leq k$) has the same parity as the entry ${n+t \choose j}$ of Pascal's triangle, by  \CountingLemmaParity. 
Now, $k=2^{q-1}+m-1$ while $n+t = 2^{q}+2m = 2(2^{q-1}+m)$. 
The central coefficient ${2(2^{q-1}+m) \choose 2^{q-1}+m}$ of the $(n+t)^{\mbox{\tiny th}}$ row of Pascal's triangle is even by \CountingLemmaCentral, so this entry does not contribute to the tally of odd numbers on this row. 
Therefore the number of odds in row $n$ of our array is exactly one-half the number of odds on the $(n+t)^{\mbox{\tiny th}}$ row of Pascal's triangle, which is $\frac{1}{2} \cdot 2^{d(n+t)}$. 
This completes the confirmation of our claim when $n$ is odd.\hfill\QED

\vspace*{0.1in}
{\bf \S 5\ \  Some thoughts on extending this work.} 
Given that our proofs employ elementary techniques, perhaps these proofs can be modified to obtain more general results (which is, indeed, a crucial function of rigorous proof in mathematics). 
One possible direction is to consider patterns in truncated Pascal arrays modulo other primes or prime powers, cf.\ \cite{Granville}. 
Also, when $t=1$, the truncated Pascal array is just the Catalan triangle. 
In this case, as mentioned in \S 3, the nonzero numbers in any given row are known to be dimensions of certain fundamental representations of the associated symplectic Lie group. 
With $t=1$, the $t$-admissible $(n,k)$-columnar tableaux we presented in \S 3 coincide (after a simple change in the alphabet of tableaux entries) with columnar symplectic tableaux of \cite{Proctor2}. 
It might be interesting to consider what similar algebraic contexts might be found for $t$-admissible $(n,k)$-columnar tableaux when $t > 1$.

\renewcommand{\refname}{\normalsize \bf References}
\renewcommand{\baselinestretch}{1}

\end{document}